\documentclass[11pt]{article}
\renewcommand{\baselinestretch}{1.2}
\newcommand{\dated}{\mbox{} \hfill {\small [{\tt \today}]}} \usepackage{amsmath,amssymb,amscd}
\newcommand{\pf}[1]{\trivlist \item[\hskip\labelsep\it #1\ ]}
\newcommand{\varpf}[1]{\trivlist \item[\hskip\labelsep\sc #1:]}
\newcommand{\qedbox}{$\rlap{$\sqcap$}\sqcup$}
\newcommand{\qed}{\qquad \qedbox \endtrivlist}
\newcommand{\varqed}{\hfill \rule{0.6em}{0.6em} \endtrivlist}
\newenvironment{proof}{\pf{Proof}}{\qed}

\newenvironment{remark}{\pf{Remark}}{\endtrivlist}
\newenvironment{remarks}{\pf{Remarks} 
   \begin{enumerate}}{\end{enumerate} \endtrivlist}
\newenvironment{example}{\pf{Example}}{\endtrivlist}
\newenvironment{examples}{\pf{Examples} 
   \begin{enumerate}}{\end{enumerate} \endtrivlist}
\newenvironment{items}{
  \begin{enumerate} 
                    
  }{\end{enumerate}}
\newenvironment{alphitems}{
  \begin{enumerate} 
                    
  }{\end{enumerate}}
\newenvironment{keywords}{\noindent\small {\it Keywords\/}:}{\vskip 4pt}
\newenvironment{classification}{\noindent\small 2000 {\it Mathematics Subject
Classification\/}:}{\vskip 12pt}

\newcommand{\comps}{{\mathbb C}}

\newcommand{\posints}{{\mathbb N}}

\newcommand{\tensor}{\otimes}
\newcommand{\ttensor}{\tilde{\otimes}}

\newcommand{\cstar}{{C^\ast}}

\newcommand{\id}{{\mathrm{id}}}

\newcommand{\A}{{\mathfrak A}}

\newcommand{\varcl}[1]{\overline{#1}}

\newtheorem{theorem}{Theorem}[section]
\newtheorem{lemma}[theorem]{Lemma}
\newtheorem{corollary}[theorem]{Corollary}
\newtheorem{proposition}[theorem]{Proposition}
\newtheorem{df}[theorem]{Definition}
\newenvironment{definition}{\begin{df} \rm}{\end{df}}

\newcommand{\QSL}{\mathit{QSL}}
\newcommand{\Rep}{\operatorname{Rep}}
\newcommand{\Cyc}{\operatorname{Cyc}}
\newcommand{\PF}{\operatorname{PF}}
\newcommand{\UPF}{\operatorname{UPF}}
\title{Representations of locally compact groups on $\QSL_p$-spaces
  and a $p$-analog of the Fourier--Stieltjes algebra}
\author{\it Volker Runde\thanks{Research supported by NSERC under grant no.\ 227043-00.}}
\date{}
\begin{document}
\maketitle
\begin{abstract}
For a locally compact group $G$ and $p \in (1,\infty)$, we define $B_p(G)$ to be the space of all coefficient functions of isometric representations of $G$ on quotients of subspaces of
$L_p$ spaces. For $p =2$, this is the usual Fourier--Stieltjes algebra. We show that $B_p(G)$ is a commutative Banach algebra that contractively (isometrically, if $G$ is amenable) contains the
Fig\`a-Talamanca--Herz algebra $A_p(G)$. If $2 \leq q \leq p$ or $p \leq q \leq 2$, we have a contractive inclusion $B_q(G) \subset B_p(G)$. We also show that $B_p(G)$ embeds contractively into the multiplier
algebra of $A_p(G)$ and is a dual space. For amenable $G$, this multiplier algebra and $B_p(G)$ are isometrically isomorphic.
\end{abstract}
\section*{Introduction}
In \cite{Eym}, P.\ Eymard introduced the {\it Fourier algebra\/} $A(G)$ of a locally
compact group $G$. If $G$ is abelian with dual group $\Gamma$, the Fourier
transform yields an isometric isomorphism of $L_1(\Gamma)$ and $A(G)$: this
motivates (and justifies) the name.
\par
For any $p \in (1,\infty)$, the {\it Fig\`a-Talamanca--Herz algebra\/}
$A_p(G)$ is defined as the collection of those functions $f \!: G \to \comps$
such that there are sequences $( \xi_n)_{n=1}^\infty$ in $L_{p'}(G)$ --- $p'
\in(1,\infty)$ being such that $\frac{1}{p} + \frac{1}{p'} =1$ --- and
$( \phi_n)_{n=1}^\infty$ in $L_p(G)$ such that
\begin{equation} \label{Ap1}
  f(x) = \sum_{n=1}^\infty \langle \lambda_{p'}(x) \xi_n, \phi_n \rangle
  \qquad (x \in G),
\end{equation}
where $\lambda_{p'}$ denotes the regular left representation of $G$ on
$L_{p'}(G)$, and
\begin{equation} \label{Ap2}
  \sum_{n=1}^\infty \| \xi_n \| \| \phi_n \| < \infty.
\end{equation}
The norm of $f \in A_p(G)$ is the infimum over all expressions of the form
(\ref{Ap2}) satisfying (\ref{Ap1}). These Banach algebras were first considered by C.\
Herz (\cite{Her1} and \cite{Her2}); their study has been an active area of
research ever since (\cite{Cow}, \cite{For1}, \cite{For2}, \cite{LNR}, \cite{Mia}, and
many more). For $p =2$, the algebra $A_p(G)$ is nothing but the Fourier
algebra $A(G)$.
\par
Another algebra introduced by Eymard in \cite{Eym} is the {\it
  Fourier--Stieltjes algebra\/} $B(G)$. For abelian $G$, it is isometrically
isomorphic to $M(\Gamma)$ via the Fourier--Stieltjes transform. It consists of
all coefficient functions of unitary representations of $G$ on some Hilbert
space and contains $A(G)$ as a closed ideal.
\par
Is there, for general $p \in (1,\infty)$, an analog of $B(G)$ in a $p$-setting 
that relates to $A_p(G)$ as does $B(G)$ to $A(G)$?
\par
In the literature (see, e.g., \cite{Cow}, \cite{For2}, \cite{Mia},
\cite{Pie}), sometimes an algebra $B_p(G)$ is considered: it is defined as the
multiplier algebra of $A_p(G)$. If $p =2$ and if $G$ is amenable, we do have
$B(G) =B_p(G)$; for non-amenable $G$, however, $B(G) \subsetneq B_2(G)$
holds. Hence, the value of $B_p(G)$ as the appropriate replacement for $B(G)$
when dealing with $A_p(G)$ is {\it a priori\/} limited to the amenable case.
\par
In the present paper, we pursue a novel approach. We define $B_p(G)$ to
consist of the coefficient functions of all representations of $G$ on
quotients of subspaces of $L_{p'}$-spaces, so-called $\QSL_{p'}$-spaces. This
class of spaces is identical with the $p'$-spaces considered in \cite{Her2}
and turns out to be appropriate for our purpose (such representations were
considered only recently, in a completely different context, in \cite{JM}).
\par
We list some properties of our $B_p(G)$:
\begin{itemize}
\item Under pointwise multiplication, $B_p(G)$ is a commutative Banach algebra
  with identity.
\item $A_p(G)$ is an ideal of $B_p(G)$, into which it contractively embeds
  (isometrically if $G$ is amenable).
\item If $2 \leq q \leq p$ or $p \leq q \leq 2$, we have a contractive
  inclusion of $B_q(G)$ in $B_p(G)$.
\item $B_p(G)$ is a dual Banach space.
\item $B_p(G)$ embeds contractively into the multiplier algebra of $A_p(G)$
and is isometrically isomorphic to it if $G$ is amenable. 
\end{itemize} 
\par
This list shows that our $B_p(G)$ relates to $A_p(G)$ in a fashion similar to
how $B(G)$ relates to $A(G)$ and therefore may be the right substitute for
$B(G)$ when working with Fig\`a-Talamanca--Herz algebras. 
\par
The main challenge when defining $B_p(G)$ and trying to establish its
properties is that the powerful methods from $\cstar$- and von Neumann
algebras are no longer at one's disposal for $p \neq 2$, so that one has to look for
appropriate substitutes.
\section{Group representations and $\QSL_p$-spaces}
We begin with defining what we mean by a representation of a locally compact group on a Banach space:
\begin{definition} \label{repdef}
A {\it representation\/} of a locally compact group $G$ (on a Banach space) is a pair
$(\pi,E)$ where $E$ is a Banach space and $\pi$ is a group homomorphism from $G$ into the invertible isometries on $E$ which is continuous with respect to
the given topology on $G$ and the strong operator topology on ${\cal B}(E)$.
\end{definition}
\begin{remarks}
\item Our definition is more restrictive than the usual definition of a representation, which does not require the range of $\pi$ to consist of isometries. Since we will not encounter any other representations, however,
we feel justified to use the general term ``representation'' in the sense defined in Definition \ref{repdef}.
\item Any representation $(\pi,E)$ of a locally compact group $G$ induces a representation of the group algebra $L^1(G)$ on $E$, i.e.\ a contractive algebra homomorphism $L_1(G)$ to ${\cal B}(E)$ --- which we shall denote
likewise by $\pi$ ---  through
\begin{equation} \label{integral}
  \pi(f) :=\int_G f(x) \pi(x) \, dx \qquad (f \in L^1(G)),
\end{equation}
where the integral (\ref{integral}) converges with respect to the strong operator topology. 
\item Instead of requiring $\pi$ to be continuous with respect to the strong operator topology on ${\cal B}(E)$, we could have demanded that $\pi$ be continuous with respect to the weak operator topology on
${\cal B}(E)$: both definitions are equivalent by \cite{GL}.
\end{remarks}
\begin{definition} \label{repdefs}
Let $G$ be a locally compact group, and let $(\pi, E)$ and $(\rho, F)$ be representations of $G$. Then:
\begin{alphitems}
\item $(\pi,E)$ and $(\rho,F)$ are said to be {\it equivalent\/} if there is an invertible isometry $V \!: E \to F$ such that
\[
  V \pi(x) V^{-1} = \rho(x) \qquad (x \in G).
\]
\item $(\rho,F)$ is called a subrepresentation of $(\pi,E)$ if $F$ is a closed subspace of $E$ such that
\[
  \rho(x) = \pi(x) |_F \qquad (x \in G).
\]
\item $(\rho,F)$ is said to be {\it contained\/} in $(\pi,E)$ --- in symbols: $(\rho,F) \subset (\pi,E)$ --- if $(\rho,F)$ is equivalent to a subrepresentation of $(\pi,E)$.
\end{alphitems}
\end{definition}
\par
Throughout, we shall often not tell a particular representation apart from its equivalence class. This should, however, not be a source of confusions.
\par
In this paper, we are interested in representations of locally compact groups on rather particular Banach spaces:
\begin{definition} \label{sql}
Let $p \in (1,\infty)$.
\begin{alphitems}
\item A Banach space is called an $L_p$-space if it is of the form $L_p(X)$ for some measure space $X$.
\item A Banach space is called a $\QSL_p$-space if it is isometrically isomorphic to a quotient of a subspace of an $L_p$-space.
\end{alphitems}
\end{definition}
\begin{remarks}
\item Equivalently, a Banach space is a $\QSL_p$-space if and only if it is a subspace of a quotient of an $L_p$-space.
\item Trivially, the class of $\QSL_p$-spaces is closed under taking subspaces and quotients.  
\item If $(E_\alpha )_\alpha$ is a family of $\QSL_p$-spaces, its $\ell_p$-direct sum $\text{$\ell_p$-}\bigoplus_\alpha E_\alpha$ is again a $\QSL_p$-space.
\item If $E$ is a $\QSL_p$-space and if $p' \in (1,\infty)$ is such that $\frac{1}{p} + \frac{1}{p'} =1$, the dual space $E^\ast$ is an $\QSL_{p'}$-space. In particular, every $\QSL_p$-space is reflexive.
\item By \cite[\S 4, Theorem 2]{Kwa}, the $QSL_p$-spaces are precisely the $p$-spaces in the sense of \cite{Her1}, i.e.\ those Banach spaces $E$ such that for any two measure spaces $X$ and $Y$ the amplification
map
\[
  {\cal B}(L_p(X),L_p(Y)) \to {\cal B}(L_p(X,E),L_p(Y,E)), \quad T \mapsto T \tensor \id_E
\]
is an isometry. In particular, an $L_q$-space is a $\QSL_p$-space if and only if
$2 \leq q \leq p$ or $p \leq q \leq 2$. Consequently, if $2 \leq q \leq p$ or $p \leq q \leq 2$, then every $\QSL_q$-space is a $\QSL_p$-space.
\item All ${\cal L}_p$-spaces in the sense of \cite{LR} --- and, more generally, all $\mathfrak{L}^g_p$-spaces in the sense of \cite{DF} --- are $\QSL_p$-spaces.
\item Since the class of $L_p$-space is stable under forming ultrapowers (\cite{Hei}), so is the class of $\QSL_p$-spaces (this immediately yields that $\QSL_p$-spaces are not only reflexive, but actually superreflexive). 
In the case where $X = Y= \comps$, the $\QSL_p$-spaces are therefore precisely those that occur in \cite[Theorem 4.1]{LeM} and play the r\^ole played by Hilbert spaces in Ruan's representation theorem for operator spaces (\cite[Theorem 2.3.5]{ER}).
\end{remarks}
\section{The linear space $B_p(G)$}
We shall not so much be concerned with representations themselves, but rather with certain functions associated with them:
\begin{definition}
Let $G$ be a locally compact group, and let $(\pi,E)$ be a representation of $G$. A {\it coefficient function\/} of $(\pi,E)$ is a function $f \!: G \to \comps$ of the form
\begin{equation} \label{coeff}
  f(x) = \langle \pi(x) \xi, \phi \rangle \qquad (x \in G),
\end{equation}
where $\xi \in E$ and $\phi \in E^\ast$.
\end{definition}
\begin{remark}
It is clear that every coefficient function of the form (\ref{coeff}) must be both bounded --- by $\| \xi \| \| \phi \|$ --- and continuous. 
\end{remark}
\par
For any locally compact group $G$ and $p \in (1,\infty)$, we denote by $\Rep_p(G)$ the collection of all (equivalence classes) of representations of $G$ on a $\QSL_p$-space.
\begin{examples}
\item The {\it left regular representation\/} $(\lambda_p,L_p(G))$ of $G$ with
\[
  \lambda_p(x)\xi(y) := \xi(x^{-1}y) \qquad (x,y \in G, \, \xi \in L_p(G))
\]
belongs to $\Rep_p(G)$.
\item For any $QSL_p$-space $E$, the {\it trivial representation\/} $(\id_E,E)$ lies in $\Rep_p(G)$.
\item For $2 \leq q \leq p$ or $p \leq q \leq 2$, we have $\Rep_q(G) \subset \Rep_p(G)$, so that, in particular, every unitary representation of $G$ on a Hilbert space belongs to $\Rep_p(G)$.
\end{examples}
\par
We can now define the main object of study in this article:
\begin{definition} \label{Bpdef}
Let $G$ be a locally compact, let $p \in (1,\infty)$, and let $p' \in (1,\infty)$ be such that $\frac{1}{p} + \frac{1}{p'} =1$. Let
\[
  B_p(G) := \left\{ f \!: G \to \comps : \text{$f$ is a coefficient of some $(\pi,E) \in \Rep_{p'}(G)$} \right\}.
\]
\end{definition}
\begin{remarks}
\item In the literature (see, for instance, \cite{Pie}), the symbol $B_p(G)$ is usually used to denote the {\it multiplier algebra\/} of $A_p(G)$, i.e.\ the set of those continuous functions $f$ on $G$ such that
$f A_p(G) \subset A_p(G)$.
\item Since subspaces and quotients of Hilbert spaces are again Hilbert spaces, $B_2(G)$ is just the usual Fourier--Stieltjes algebra $B(G)$ introduced in \cite{Eym}. For amenable $G$, this is consistent 
with the usage in \cite{Pie}. In the non-amenable case, however, $B_2(G) = B(G)$ as defined in Definition \ref{Bpdef} and $B_2(G)$ in the sense of \cite{Pie} denote different objects.
\end{remarks}
\par
We conclude this section with proving a few, very basic properties of $B_p(G)$:
\begin{lemma} \label{Bplem}
Let $G$ be a locally compact group,  let $p \in (1,\infty)$, let $p' \in (1,\infty)$ be such that $\frac{1}{p} + \frac{1}{p'} =1$, and let $f \!:G \to \comps$ be a function such that the following holds: There are sequences
$((\pi_n,E_n) )_{n=1}^\infty$, $( \xi_n )_{n=1}^\infty$, and $( \phi_n )_{n=1}^\infty$ with $(\pi_n,E_n) \in \Rep_{p'}(G)$, $\xi_n \in E_n$, and $\phi_n \in E_n^\ast$ for $n \in \posints$ such that
\[
  \sum_{n=1}^\infty \| \xi_n \| \| \phi_n \| < \infty
\]
and
\[
  f(x) = \sum_{n=1}^\infty \langle \pi_n(x) \xi_n, \phi_n \rangle \qquad (x \in G).
\]
Then $f$ lies in $B_p(G)$.
\end{lemma}
\begin{proof}
Without loss of generality, we may suppose that
\[
  \sum_{n=1}^\infty \| \xi_n \|^{p'} < \infty \qquad\text{and}\qquad  \sum_{n=1}^\infty \| \phi_n \|^p < \infty.
\]
\par
Define $(\pi,E) \in \Rep_{p'}(G)$ by letting $E := \text{$\ell_{p'}$-}\bigoplus_{n=1}^\infty E_n$ and, for $\eta = (\eta_1, \eta_2, \ldots) \in E$,
\[
  \pi(x)\eta := (\pi_1(x)\eta_1, \pi_2(x),\eta_2, \ldots) \qquad (x \in G).
\]
It follows that $\xi := (\xi_1, \xi_2, \ldots) \in E$, that $\phi := (\phi_1, \phi_2, \ldots) \in E^\ast$, and that $f$ is a coefficient function of $(\pi,E)$ --- therefore belonging to $B_p(G)$.
\end{proof}
\par
For any topological space $\Omega$, we use ${\cal C}_\mathrm{b}(\Omega)$ to denote the bounded continuous functions on it.
\begin{proposition} \label{Bplin}
Let $G$ be a locally compact group, and let $p \in (1,\infty)$. Then $B_p(G)$ is a linear subspace of ${\cal C}_\mathrm{b}(G)$ containing $A_p(G)$. Moreover, if $2 \leq q \leq p$ or $p \leq q \leq 2$, we have $B_q(G) \subset B_p(G)$.
\end{proposition}
\begin{proof}
We have already seen that $B_p(G) \subset {\cal C}_\mathrm{b}(G)$.
\par
Let $p' \in (1,\infty)$ be such that $\frac{1}{p} + \frac{1}{p'} =1$, and let $f_1, f_2 \in B_p(G)$. By the definition of $B_p(G)$, there are $(\pi_1, E_1), (\pi_2,E_2) \in \Rep_{p'}(G)$ such that $f_j$ is a coefficient function of 
$(\pi_j, E_j)$  for $j=1,2$. It is clear that the pointwise sum $f_1 + f_2$ is then of the form considered in Lemma \ref{Bplem} (take $\xi_3 = \xi_4 = \cdots = 0$) and thus contained in $B_p(G)$.
\par
To see that $A_p(G) \subset B_p(G)$, apply Lemma \ref{Bplem} again with $(\pi_n,E_n) = (\lambda_{p'}, L_{p'}(G))$ for $n \in \posints$.
\par
Let $2 \leq q \leq p$ or $p \leq q \leq 2$, and let $q' \in (1,\infty)$ be such that $\frac{1}{q} + \frac{1}{q'} = 1$. Since every $\QSL_{q'}$ space is a is a $\QSL_{p'}$-space, the inclusion $B_q(G) \subset B_p(G)$ holds.
\end{proof}
\section{Tensor products of $\QSL_p$-spaces}
Let $G$ be a locally compact group. In $B(G) = B_2(G)$, the pointwise product of functions corresponds to the tensor product of representations, which, in turn, relies on the existence of the Hilbert space tensor product. 
In order to turn $B_p(G)$ into an algebra, we will therefore equip, in this section, the algebraic tensor product of two $\QSL_{p'}$-spaces, where $\frac{1}{p} + \frac{1}{p'} =1$, with a suitable norm.
\par
The main result is the following:
\begin{theorem} \label{tensorthm}
Let $E$ and $F$ be $\QSL_p$-spaces. Then there is a norm $\| \cdot \|_p$ on the algebraic tensor product $E \tensor F$ such that:
\begin{items}
\item $\| \cdot \|_p$ dominates the injective norm;
\item $\| \cdot \|_p$ is a cross norm;
\item the completion $E \ttensor_p F$ of $E \tensor F$ with respect to $\| \cdot \|_p$ is a $\QSL_p$-space.
\end{items}
Moreover, if $G$ is a locally compact group with $(\pi,E),(\rho,F) \in \Rep_p(G)$, then $(\pi \tensor \rho, E \ttensor_p F) \in \Rep_p(G)$ is well defined through
\[
  (\pi(x) \tensor \rho(x))(\xi \tensor \eta) := \pi(x) \xi \tensor \rho(x)\eta \qquad (x \in G, \, \xi \in E, \, \eta \in F).
\]
\end{theorem}
\begin{proof}
Let $X$ be a measure space, let $E_1$ and $F_1$ be closed subspaces of $L_p(X)$, and let $E_2$ and $F_2$ be closed subspaces of $E_1$ and $F_1$, respectively, such that $E = E_1/E_2$ and $F = F_1 / F_2$.
\par
We may embed the algebraic tensor product $L_p(X) \tensor L_p(X)$ into the vector valued $L_p$-space $L_p(X,L_p(X))$ and thus equip it with a norm, which we denote by $\|| \cdot \||_p$ which dominates the injective norm
on $L_p(X) \tensor L_p(X)$ (\cite[7.1, Proposition]{DF}). Of course, we may restrict $\|| \cdot \||_p$ to $E_1 \tensor E_2$. We denote the (uncompleted) injective tensor product by $\tensor_\epsilon$. Since $\tensor_\epsilon$
respects passage to subspaces, we see that the identity on $E_1 \tensor F_1$ induces a contraction from $(E_1 \tensor E_2, \|| \cdot \||_p)$ to $E_1 \tensor_\epsilon F_1$. Let $\pi_E \!: E_1 \to E$ and $\pi_F \!: F_1 \to F$ denote
the canonical quotient maps. The mapping property of the injective tensor product then yields that
\[
  \pi_E \tensor \pi_F \!: (E_1 \tensor F_1, \|| \cdot \||_p) \to E_1 \tensor_\epsilon F_1 \to E \tensor_\epsilon F
\]
is a surjective contraction, so that, in particular, $\ker(\pi_E \tensor \pi_F)$ is closed in $(E_1 \tensor F_1, \|| \cdot \||_p)$. Let $\| \cdot \|_p$ denote the induced quotient norm on $E \tensor F = (E_1 \tensor F_1) /
\ker(\pi_E \tensor \pi_F)$. It is immediate that $\| \cdot \|_p$ dominates the injective tensor norm on $E \tensor F$, so that (i) holds. Moreover, since $\|| \cdot \||_p$ is a cross norm on $E_1 \tensor E_2$, it is clear that
$\| \cdot \|_p$ is at least subcross on $E \tensor F$. Since $\| \cdot \|_p$, however, dominates the injective norm --- which is a cross norm --- on $E \tensor F$, we conclude that $\| \cdot \|_p$ is indeed a cross norm on 
$E \tensor F$. This proves (ii).
\par
For notational convenience, we write $L_p(X) \tensor_p L_p(X) := (L_p(X) \tensor L_p(X), \|| \cdot \||_p)$, and let $E \tensor_p F := (E \tensor F, \| \cdot \|_p)$. Let $Y$ and $Z$ be any measure spaces. In view of 
\cite[7.2 and 7.3]{DF}, it is clear that the amplification map
\[
  {\cal B}(L_p(Y),L_p(Z)) \to {\cal B}(L_p(Y,L_p(X) \tensor_p L_p(X)),L_p(Z,L_p(X) \tensor_p L_p(X))), \quad T \to T \tensor \id
\]
is an isometry, and from \cite[7.4, Proposition]{DF}, we conclude that the same is true for
\begin{equation} \label{pspace}
  {\cal B}(L_p(Y),L_p(Z)) \to {\cal B}(L_p(Y,E \tensor_p F),L_p(Z,E \tensor_p F)), \quad T \to T \tensor \id.
\end{equation}
However, if we replace $E \tensor_p F$ in (\ref{pspace}) by its completion $E \ttensor_p F$, (\ref{pspace}) obviously remains an isometry. Hence, $E \ttensor_p F$ is a $p$-space in the sense of \cite{Her1} and thus a
$\QSL_p$-space by \cite[\S4, Theorem 2]{Kwa}.
\par
For the moreover part of the theorem, it is sufficient to show that, for $S \in {\cal B}(E)$ and $T \in {\cal B}(F)$, their tensor product $S \tensor T$ is continuous on $E \tensor_p F$ and has operator norm at most $\| S \| \| T \|$.
We first treat the case where $S = \id_E$. Let $E_1 \tensor_p F$ stand for $E_1 \tensor F$ equipped with the norm obtained by factoring $E_1 \tensor F_2$ out of $(E_1 \tensor F_1, \|| \cdot \||_p)$. From \cite[7.3]{DF}, it follows 
that $\id_{E_1} \tensor T \in {\cal B}(E_1 \tensor F)$ and has operator norm such that
\[
  \| \id_{E_1} \tensor T \|_{{\cal B}(E_1 \tensor_p F)} = \| T \|_{{\cal B}(F)}.
\]
It is easy to see that $E \tensor F$ is, in fact, the quotient space of $E_1 \tensor_p F$ module $E_2 \tensor F$, it follows that
\[
  \| \id_E \tensor T \|_{{\cal B}(E \tensor_p F)} \leq \| \id_{E_1} \tensor_p T \|_{{\cal B}(E_1 \tensor F)} = \| T \|_{{\cal B}(F)}.
\]
By symmetry, we obtain that
\[
  \| S \tensor \id_F \|_{{\cal B}(E \tensor_p F)} \leq \| S \|_{{\cal B}(E)}
\]
as well. Consequently,
\[
  \| S \tensor T \|_{{\cal B}(E \tensor_p F)} \leq \| S \tensor \id_F \|_{{\cal B}(E \tensor_p F)} \| \id_E \tensor T \|_{{\cal B}(E \tensor_p F)} \leq \| S \|_{{\cal B}(E)} \| T \|_{{\cal B}(F)}
\]
holds.
\end{proof}
\begin{remarks}
\item For a measure space $X$ and for a $\QSL_p$-space $E$, the tensor product $L_p(X) \ttensor_p E$ constructed in the proof of Theorem \ref{tensorthm} is nothing but the vector valued $L_p$-space $L_p(X,E)$.
\item We suspect, but have been unable to prove, that $\| \cdot \|_p$ is the Chevet--Saphar tensor norm $d_p$ on $E \tensor F$ (see \cite[12.7]{DF}). This is indeed the case when both $E$ and $F$ are 
$\mathfrak{L}_p^g$-spaces (see \cite{JM}).
\end{remarks}
\par
We conclude this section with two simple corollaries of Theorem \ref{tensorthm}:
\begin{corollary} \label{tensorcor1}
Let $G$ be a locally compact group, let $p \in (1,\infty)$, and let $f, g \!: G \to \comps$ be coefficient function of $(\pi,E)$ and $(\rho,F)$ in $\Rep_p(G)$, respectively, namely
\[
  f(x) = \langle \pi(x)\xi, \phi \rangle \quad\text{and}\quad g(x) = \langle \rho(x)\eta, \psi \rangle \qquad (x \in G)
\]
where $\xi \in E$, $\phi \in E^\ast$, $\eta \in F$, and $\psi \in F^\ast$. Then $\phi \tensor \psi \!: E \tensor F \to \comps$ is continuous with respect to $\| \cdot \|_p$ with norm at most $\| \phi \| \| \psi \|$, so that the 
pointwise product of $f$ and $g$ is a coefficient function of $(\pi \tensor \rho, E \ttensor_p F)$, namely
\[
  f(x) g(x) = \langle (\pi(x) \tensor \rho(x))(\xi \tensor \eta), \phi \tensor \psi \rangle \qquad (x \in G).
\]
\end{corollary}
\begin{proof}
In view of the definition of $(\pi \tensor \rho, E \ttensor_p F)$, only the claim about $\phi \tensor \psi$ needs some consideration: it is, however, an immediate consequence of Theorem \ref{tensorthm}(i) and (ii).
\end{proof}
\begin{corollary} \label{tensorcor2}
Let $G$ be a locally compact group, and let $p \in (1,\infty)$. Then $B_p(G)$ is a unital subalgebra of ${\cal C}_\mathrm{b}(G)$.
\end{corollary}
\begin{proof}
By Proposition \ref{Bplin}, $B_p(G)$ is a linear subspace of ${\cal C}_\mathrm{b}(G)$, and by Corollary \ref{tensorcor1}, it is a subalgebra. The constant function $1$ is a coefficient function of any trivial representation of
$G$ on an $\QSL_p$-space.
\end{proof}
\section{The Banach algebra $B_p(G)$}
Our next goal is to equip the algebra $B_p(G)$ with a norm turning it into a Banach algebra.
\begin{definition} \label{cycdef}
Let $G$ be a locally compact group, and let $(\pi,E)$ be a representation of $G$. Then $(\pi,E)$ is called {\it cyclic\/} if there is $x \in E$ such that $\pi(L_1(G))x$ is dense in $E$. For $p \in (1,\infty)$, we let
\[
  \Cyc_p(G) := \{ (\pi,E) : \text{$(\pi,E)$ is cyclic} \}.
\]
\end{definition}
\begin{remark}
Let $f \in B_p(G)$ be a coefficient function of $(\pi,E) \in \Rep_p(G)$, i.e.\
\[
  f(x) = \langle f(x)\xi, \phi \rangle \qquad (x \in G)
\]
with $\xi \in E$ and $\phi \in E^\ast$. Let $F := \varcl{\pi(L_1(G))\xi}$, and define $\rho \!: G \to {\cal B}(F)$ by restriction of $\pi(x)$ to $F$ for each $x \in G$. Then $(\rho,F)$ is cyclic with $f$ as a coefficient function.
\end{remark}
\begin{definition} \label{normdef}
Let $G$ be a locally compact group, let $p, p' \in (1,\infty)$ be dual to each other --- meaning: $\frac{1}{p} + \frac{1}{p'} = 1$ ---, and let $f \in B_p(G)$. We define $\| f \|_{B_p(G)}$ as the infimum over all expressions
$\sum_{n=1}^\infty \| \xi_n \| \| \phi_n \|$, where, for each $n \in \posints$, there is $(\pi_n,E_n) \in \Cyc_{p'}(E)$ with $\xi_n \in E_n$ and $\phi_n \in E_n^\ast$ such that
\[
  \sum_{n=1}^\infty \| \xi_n \| \| \phi_n \| < \infty \qquad\text{and}\qquad 
  f(x) = \sum_{n=1}^\infty \langle \pi_n(x) \xi_n, \phi_n \rangle \quad (x \in G).
\]
\end{definition}
\begin{remarks}
\item In view of the remark after Definition \ref{cycdef}, it is clear that $\| \cdot \|_{B_p(G)}$ is well defined, and it is easily checked that $\| \cdot \|_{B_p(G)}$ is indeed a norm on $B_p(G)$.
\item One might think that it would be more appropriate to define $\| \cdot \|_{B_p(G)}$ in such a way that the infimum is taken over general $(\pi_n,E_n) \in \Rep_{p'}(G)$ instead of only in $\Cyc_{p'}(G)$. The problem here, however,
is that $\QSL_p$-spaces can be of arbitrarily large cardinality, so that $\Rep_{p'}(G)$ is not a set, but only a class. Since, for $(\pi,E) \in \Cyc_{p'}(G)$, the space $E$ has a cardinality not larger than $|L_1(G)|^{\aleph_0}$, it
follows that $\Cyc_{p'}(G)$ --- unlike all of $\Rep_{p'}(G)$ --- is indeed a set, so that it makes sense to take an infimum over it.
\end{remarks}
\par
In view of the last one of the two preceding remarks, the following lemma is comforting:
\begin{lemma} \label{cyclem1}
Let $G$ be a locally compact group, let $p,p' \in (1,\infty)$ be dual to each other, and let $(( \pi_n,E_n) )_{n=1}^\infty$ be a sequence in $\Rep_{p'}(G)$ such that, with $\xi_n \in E_n$ and $\phi_n \in E_n^\ast$ for $n \in \posints$, we have
$\sum_{n=1}^\infty \| \xi_n \| \| \phi_n \| < \infty$. Then, for each $n \in \posints$, there are $(\rho_n,F_n) \in \Cyc_{p'}(G)$ with $(\rho_n,F_n) \subset (\pi_n,E_n)$, $\eta_n \in F_n$, and $\psi_n \in E^\ast$, such that
\[
  \sum_{n=1}^\infty \| \eta_n \| \| \psi_n \| \leq \sum_{n=1}^\infty \| \xi_n \| \| \phi \|
\]
and
\[
  \sum_{n=1}^\infty \langle \rho_n(x) \eta_n, \psi_n \rangle = \sum_{n=1}^\infty \langle \rho_n(x) \xi_n, \phi_n \rangle \quad (x \in G)
\]
\end{lemma}
\begin{proof}
We proceed as in the remark immediately following Definition \ref{cycdef}: For $n \in \posints$, let $F_n := \varcl{\pi_n(L_1(G))\xi_n}$, define $\rho_n$ through restriction, let $\eta_n := \xi_n$, and let $\psi_n$ be the restriction of
$\phi_n$ to $F_n$.
\end{proof}
\begin{lemma} \label{cyclem2}
Let $G$ be a locally compact group, let $p,p' \in (1,\infty)$ be dual to each other, and let $f \in A_p(G)$. Then $\| f \|_{A_p(G)}$ is the infimum over all expressions
$\sum_{n=1}^\infty \| \xi_n \| \| \phi_n \|$, where, for each $n \in \posints$, there is $(\pi_n,E_n) \in \Cyc_{p'}(E)$ contained in $(\lambda_{p'}, L_{p'}(G))$ with $\xi_n \in E_n$ and $\phi_n \in E_n^\ast$ such that
\[
  \sum_{n=1}^\infty \| \xi_n \| \| \phi_n \| < \infty \qquad\text{and}\qquad 
  f(x) = \sum_{n=1}^\infty \langle \pi_n(x) \xi_n, \phi_n \rangle \quad (x \in G).
\]
\end{lemma}
\begin{proof}
From Lemma \ref{cyclem1}, it follows that the infimum in the statement of Lemma \ref{cyclem2} is less or equal to $\| f \|_{A_p(G)}$. Let this infimum be denoted by $C_f$. Let $\epsilon > 0$, and choose a sequence $((\pi_n,E_n))_{n=1}^\infty$ of
cyclic subrepresentations of $(\lambda_{p'}, L_{p'}(G))$ and, for each $n \in \posints$, $\xi_n \in E_n$ and $\phi_n \in E_n^\ast$ such that
\[
  \sum_{n=1}^\infty \| \xi_n \| \| \phi_n \| < C_f + \epsilon \qquad\text{and}\qquad 
  f(x) = \sum_{n=1}^\infty \langle \pi_n(x) \xi_n, \phi_n \rangle \quad (x \in G).
\]
For each $n \in \posints$, use the Hahn--Banach theorem to extend $\phi_n \in E_n^\ast$ to $\psi_n \in L_{p'}(G)^\ast = L_p(G)$ with $\| \psi_n \| = \| \phi_n \|$. It follows that
\[
  \| f \|_{A_p(G)} \leq \sum_{n=1}^\infty \| \xi_n \| \| \psi_n \| = \sum_{n=1}^\infty \| \xi_n \| \| \phi_n \| < C_f + \epsilon.
\]
Since $\epsilon > 0$ was arbitrary, we conclude that $\| f \|_{A_p(G)} \leq C_f$.
\end{proof}
\begin{definition} \label{unidef}
Let $G$ be a locally compact group, and let $p \in (1,\infty)$. Then
$(\pi,E) \in \Rep_p(G)$ is called {\it $p$-universal\/} if $(\rho,F)
\subset (\pi,E)$ for all $(\rho,F) \in \Cyc_p(G)$.
\end{definition}
\begin{example}
Let $G$ be a locally compact group, and let $p \in (1,\infty)$. Since
$\Cyc_p(G)$ is a set, we can form the $\ell_p$-direct sum of all
$(\rho,F) \in \Cyc_p(G)$. This representation is then obviously
$p$-universal.
\end{example}
\begin{lemma} \label{cyclem3}
Let $G$ be a locally compact group, let $p, p' \in (1,\infty)$ be dual
to each other, and let $(\pi,E) \in \Rep_{p'}(G)$ be
$p'$-universal. Then, for each $f \in B_p(G)$, the norm $\| f \|_{B_p(G)}$ is the infimum over all expressions
$\sum_{n=1}^\infty \| \xi_n \| \| \phi_n \|$ with $\xi_n \in E$ and $\phi_n \in E^\ast$ for each $n \in \posints$ such that
\[
  \sum_{n=1}^\infty \| \xi_n \| \| \phi_n \| < \infty \qquad\text{and}\qquad 
  f(x) = \sum_{n=1}^\infty \langle \pi(x) \xi_n, \phi_n \rangle \quad (x \in G).
\]
\end{lemma}
\begin{proof}
Obvious in the light of Definition \ref{unidef}.
\end{proof}
\par
In the end, we obtain:
\begin{theorem} \label{bathm}
Let $G$ be a locally compact group, let $p \in (1,\infty)$, and let $B_p(G)$ be equipped with $\| \cdot \|_{B_p(G)}$. Then:
\begin{items}
\item $B_p(G)$ is a commutative Banach algebra.
\item The inclusion $A_p(G) \subset B_p(G)$ is a contraction.
\item For $2 \leq q \leq p$ or $p \leq q \leq 2$, the inclusion $B_q(G) \subset B_p(G)$ is a contraction.
\end{items}
\end{theorem}
\begin{proof}
Let $p' \in (1,\infty)$ be dual to $p$, and let $(\pi,E) \in \Rep_{p'}(G)$ be $p'$-universal. It follows that $B_p(G)$ is a quotient space of the complete projective tensor product $E \ttensor_\pi E^\ast$ and thus complete. 
By Corollary \ref{tensorcor2}, $B_p(G)$ is an algebra, so that all that remains to prove (i) is to show that $\| \cdot \|_{B_p(G)}$ is submultiplicative.
\par
Let $f,g \in B_p(G)$, and let $\epsilon > 0$. Let $((\pi_n,E_n))_{n=1}^\infty$ and $((\rho_n,F_n))_{n=1}^\infty$ be sequences in $\Cyc_{p'}(G)$ and, for $n \in \posints$, let $\xi_n \in E_n$, $\phi_n \in E_n^\ast$, $\eta_n \in F_n$,
and $\psi_n \in F_n^\ast$ such that
\[
  f(x) = \sum_{n=1}^\infty \langle \pi_n(x) \xi_n, \phi_n \rangle \quad\text{and}\quad
  g(x) = \sum_{n=1}^\infty \langle \rho_n(x) \eta_n, \psi_n \rangle\qquad (x \in G)
\]
and
\[
  \sum_{n=1}^\infty \| \xi_n \| \| \phi_n \| \leq \| f \|_{B_p(G)} + \epsilon \qquad\text{and}\qquad \sum_{n=1}^\infty \| \eta_n \| \| \psi_n \| \leq \| g \|_{B_p(G)} + \epsilon.
\]
By the ``moreover'' part of Theorem \ref{tensorthm}, we see that $( \pi_n \tensor \rho_m, E_n \ttensor_p F_m) \in \Rep_{p'}(G)$ for $n,m \in \posints$, and Corollary \ref{tensorcor1} yields that
\[
  f(x)g(x) = \sum_{n,m=1}^\infty \langle (\pi_n(x) \tensor \rho_m(x))(\xi_n \tensor \eta_m), \phi_n \tensor \psi_m \rangle\qquad (x \in G)
\]
and that
\begin{eqnarray*}
  \sum_{n,m=1}^\infty \| \xi_n \tensor \eta_m \|_{E_n \ttensor_p F_n} \| \phi_n \tensor \psi_m \|_{(E_n \ttensor_p F_n)^\ast}
  & \leq  & \sum_{n,m=1}^\infty \| \xi_n \| \| \eta_m \| \| \phi_n \| \psi_m \| \\
  & \leq & \left( \sum_{n=1}^\infty \| \xi_n \| \| \phi_n \| \right) \left(  \sum_{m=1}^\infty \| \eta_m \| \| \psi_m \|\right) \\
  & \leq & (\| f \|_{B_p(G)} + \epsilon)(\| g \|_{B_p(G)} + \epsilon).
\end{eqnarray*}
From Lemma \ref{cyclem1} and Definition \ref{normdef}, we conclude that
\[
  \| f g \|_{B_p(G)} \leq (\| f \|_{B_p(G)} + \epsilon)(\| g \|_{B_p(G)} + \epsilon).
\]
Since $\epsilon > 0$ was arbitrary, this yields the submultiplicativity of $\| \cdot \|_{B_p(G)}$ and thus completes the proof of (i).
\par
From Lemma \ref{cyclem2} and Definition \ref{normdef}, (ii) is immediate.
\par
Let $2 \leq q \leq p$ or $p \leq q \leq 2$, and let $q' \in (1,\infty)$ be dual to $q$. Since $\Cyc_{q'}(G) \subset \Cyc_{p'}(G)$, this proves (iii).
\end{proof}
\section{$B_p(G)$ and $A_p(G)$}
For any locally compact group $G$, the Fourier algebra $A(G)$ embeds
isometrically into $B(G)$ and can be identified with the closed ideal of
$B(G)$ generated by the functions in $B(G)$ with compact support
(\cite{Eym}). 
\par
For general $p \in (1,\infty)$, the only information we have so far about the 
relation between $B_p(G)$ and $A_p(G)$ is Theorem \ref{bathm}(ii). In the 
present section, we further explore the relation between those algebras. 
\par
Our first result is known for $p=2$ as {\it Fell's absorption
  principle\/}:
\begin{proposition} \label{Fell}
Let $G$ be a locally compact group, let $p \in (1,\infty)$, and let
$(\pi,E) \in \Rep_p(G)$. Then the representations $(\lambda_p \tensor
\pi, L_p(G,E))$ and $(\lambda_p \tensor \id_E , L_p(G,E))$ are equivalent.
\end{proposition}
\begin{proof}
The proof very much goes along the lines of the case $p=2$. 
\par
Let ${\cal C}_{00}(G,E)$ denote the continuous $E$-valued functions on
$G$ with compact support (so that ${\cal C}_{00}(G,E)$ is a dense
subspace of $L_p(G,E)$). Define $W_\pi \!: {\cal C}_{00}(G,E) \to
{\cal C}_{00}(G,E)$ by letting
\[
  (W_\pi \xi)(x) := \pi(x)\xi(x) \qquad (\xi \in {\cal C}_{00}(G,E), \, x
  \in G).
\]
Since $\pi(G)$ consists of isometries, we have
\[
  \| W_\pi \xi \|_{L_p(G,E)}^p = \int_G \|  \pi(x)\xi(x) \|^p \, dx =
  \int_G \| \xi(x) \|^p \, dx
  \qquad (\xi \in {\cal C}_{00}(G,E)),
\]
so that $W_p$ is an isometry with respect to the norm of $L_p(G,E)$
and thus extends to all of $L_p(G,E)$ as an isometry. Clearly, $W_\pi$
is invertible with inverse given by
\[
  (W_\pi^{-1} \xi)(x) := \pi(x^{-1})\xi(x) \qquad (\xi \in {\cal C}_{00}(G,E), \, x
  \in G).
\]
\par
Let $\xi \in {\cal C}_{00}(G,E)$, and let $x \in G$. Then we have
\[ 
  ((\lambda_p(x) \tensor \id_E) W_\pi^{-1} \xi)(y) =
  \pi(y^{-1}x)\xi(x^{-1}y) \qquad (y \in G)
\]
and thus
\begin{eqnarray*}
  (W_\pi(\lambda_p(x) \tensor \id_E) W^{-1}_\pi \xi)(y) & = &
  \pi(y)\pi(y^{-1}x)\xi(x^{-1}y) \\
  & = & \pi(x) \xi(x^{-1}y) \\
  & = & ((\lambda_p(x) \tensor \pi(x))\xi)(y) \qquad (y \in G).
\end{eqnarray*}
Hence, 
\[
  W_\pi(\lambda_p(x) \tensor \id_E) W^{-1}_\pi = \lambda_p(x)
  \tensor \pi(x) \qquad (x \in G)
\]
holds, so that $(\lambda_p \tensor
\pi, L_p(G,E))$ and $(\lambda_p \tensor \id_E , L_p(G,E))$ are equivalent as claimed.
\end{proof}
\begin{corollary} \label{Fellcor1}
Let $G$ be a locally compact group, let $p \in G$, let $f \in A_p(G)$,
and let $g \in B_p(G)$. Then $fg$ lies in $A_p(G)$ such that
\[
  \| fg \|_{A_p(G)} \leq \| f \|_{A_p(G)} \| g \|_{B_p(G)}.
\]
\end{corollary}
\begin{proof}
Apply Proposition \ref{Fell} (with $p$ replaced by $p'$ dual to $p$)
to a $p'$-universal representation $(\pi,E) \in \Rep_{p'}(G)$. The
norm estimate is proven as is the submultiplicativity assertion of
Theorem \ref{bathm}.
\end{proof}
\par
Let $G$ be a locally compact group, and let $p \in (1,\infty)$. A {\it
  multiplier\/} of $A_p(G)$ is a function $f \in {\cal
  C}_\mathrm{b}(G)$ such that $fA_p(G) \subset A_p(G)$. We denote the
set of all multipliers of $A_p(G)$ by ${\cal M}(A_p(G))$. Clearly,
${\cal M}(A_p(G))$ is a subalgebra of ${\cal C}_\mathrm{b}(G)$. From
the closed graph theorem, it is immediate that multiplication with $f
\in {\cal M}(A_p(G))$ is a bounded linear operator on $A_p(G)$, so
that ${\cal M}(A_p(G))$ embeds canonically into ${\cal B}(A_p(G))$
turning it into a Banach algebra.
\par
We have the following (compare \cite[Lemma 0]{Her1}):
\begin{corollary} \label{Fellcor2}
Let $G$ be a locally compact group, and let $p \in (1,\infty)$. Then
$B_p(G)$ is contained in ${\cal M}(A_p(G))$ such that
\begin{equation} \label{mult1}
  \| f \|_{{\cal M}(A_p(G))} \leq \| f \|_{B_p(G)} \qquad (f \in B_p(G)).
\end{equation}
In particular, 
\begin{equation} \label{mult2}
  \| f \|_{{\cal M}(A_p(G))} \leq \| f \|_{B_p(G)}\leq \| f \|_{A_p(G)} \qquad (f \in A_p(G))
\end{equation}
holds with equality throughout if $G$ is amenable.
\end{corollary}
\begin{proof}
By Corollary \ref{Fellcor1}, $B_p(G) \subset {\cal M}(A_p(G))$ holds
as does (\ref{mult1}). The first inequality of (\ref{mult2}) follows
from (\ref{mult1}) and the second one from Theorem
\ref{bathm}(ii). Finally, if $G$ is amenable, $A_p(G)$ has an
approximate identity bounded by one (\cite[Theorem 4.10]{Pie}), so
that $\| f \|_{{\cal M}(A_p(G))} = \| f \|_{A_p(G)}$ holds for all $f
\in A_p(G)$.
\end{proof}
\begin{remark}
Let $G$ be a locally compact group such that, for any $p \in (1,\infty)$,
the embedding of $A_p(G)$ into $B_p(G)$ is an isometry. Since $A_p(G)$ is 
regular (\cite{Her2}), this means that $A_p(G)$ can be identified 
with the closed ideal of $B_p(G)$ generated by the functions in $B_p(G)$ with 
compact support. In view of Theorem \ref{bathm}(iii), this would yield a 
contractive inclusion $A_p(G) \subset A_q(G)$ whenever $2 \leq q \leq p$ or 
$p \leq q \leq 2$. Such in inclusion result is indeed true for amenable $G$ 
by C.\ Herz (\cite{Her1}) --- and also for for certain non-amenable $G$ 
(see \cite{HR}) ---, but is false for non-compact, semisimple Lie groups with 
finite center (\cite{Loh}) as was pointed out to me by Michael Cowling.
\end{remark}
\section{$B_p(G)$ as a dual space}
The Fourier--Stieltjes algebra $B(G)$ of a locally compact group $G$ can be
identified with the dual space of the full group $C^\ast$-algebra $C^\ast(G)$
(\cite{Eym}).
\par
In this section, we show that $B_p(G)$ is a dual space in a canonical fashion
for {\it arbitrary\/} $p \in (1,\infty)$. This, in turn, will enable us to
further clarify the relation between $B_p(G)$ and ${\cal M}(A_p(G))$.
\par
We begin with some more definitions:
\begin{definition} \label{pseudos}
Let $G$ be a locally compact group, let $p \in (1,\infty)$, and let
$(\pi, E) \in \Rep_p(G)$. Then:
\begin{alphitems}
\item $\| \cdot \|_\pi$ is the algebra seminorm on $L_1(G)$ defined
  through
\[
  \| f \|_\pi :=\| \pi(f) \|_{{\cal B}(E)} \qquad (f \in L_1(G)).
\]
\item The algebra $\PF_{p,\pi}(G)$ of {\it $p$-pseudofunctions
  associated with $(\pi,E)$\/} is the closure of $\pi(L_1(G))$ in ${\cal
  B}(E)$.
\item If $(\pi,E) = (\lambda_p, L_p(G))$, we simply speak of {\it
  $p$-pseudofunctions\/} and write $\PF_p(G)$ instead of  
$\PF_{p,\lambda_p}(G)$.
\item If $(\pi,E)$ is $p$-universal, we denote $\PF_{p,\pi}(G)$
  by $\UPF_p(G)$ and call it the algebra of {\it universal $p$-pseudofunctions\/}.
\end{alphitems}
\end{definition}
\begin{remarks}
\item The notion of $p$-pseudofunctions is well established in the
  literature; the other definitions seem to be new.
\item For $p=2$, the algebra $\PF_p(G)$ is the reduced group
  $\cstar$-algebra and $\UPF_p(G)$ is the full group $\cstar$-algebra
  of $G$.
\item If $(\rho,F) \in \Rep_p(G)$ is such that $(\pi,E)$ contains
  every cyclic subrepresentation of $(\rho,F)$, then $\| \cdot \|_\rho
  \leq \| \cdot \|_\pi$ holds. In particular, the definition of
  $\UPF_p(G)$ is independent of a particular $p$-universal
  representation.
\item With $\langle \cdot, \cdot \rangle$ denoting the
  $L_1(G)$-$L_\infty(G)$ duality and with $(\pi,E)$ a $p$-universal
  representation of $G$, we have
\[
  \| f \|_\pi = \sup \{ | \langle f, g \rangle | : f \in B_{p'}(G), \,
  \| g \|_{B_{p'}(G)} \leq 1 \} \qquad (f \in L_1(G)),
\]
where $p' \in (1,\infty)$ is dual to $p$: this follows from Lemma \ref{cyclem3}.
\end{remarks}
\par
We now turn to representations of Banach algebras.
\begin{definition}
A {\it representation\/} of a Banach algebra $\A$ is a pair $(\pi,E)$ where
$E$ is a Banach space and $\pi$ is a contractive algebra homomorphism from $\A$
to ${\cal B}(E)$. We call $(\pi,E)$ {\it isometric\/} if $\pi$ is an
isometry and {\it essential\/} if the linear span of $\{ \pi(a)\xi : a
\in \A, \, \xi \in E \}$ is dense in $E$.
\end{definition}
\begin{remarks}
\item As with Definition \ref{repdef}, our definition of a
  representation of a Banach algebra is somewhat more restrictive than
  the one usually used in a literature. Our reasons for this are the
  same as given after Definition \ref{repdef}.
\item If $G$ is a locally compact group and $(\pi,E)$ is a
  representation of $G$ in the sense of Definition \ref{repdef}, then
  (\ref{integral}) induces an essential representation of
  $L_1(G)$. Conversely, every essential representation of $L_1(G)$
  arises in the fashion. 
\item The notions introduced in Definition \ref{repdefs} for
  representations of locally compact groups carry over to
  representations of Banach algebras accordingly.
\end{remarks}
\par
We require three lemmas:
\begin{lemma} \label{subrep}
Let $\A$ be a Banach algebra with an approximate identity bounded by
one, and let $(\pi,E)$ be a representation of $\A$. Let $F$ be the
closed linear span of $\{ \pi(a) \xi : a \in \A, \, \xi \in E \}$, and
define
\[ 
  \rho \!: \A \to {\cal B}(F), \quad a \mapsto \pi(a) |_F.
\]
Then $(\rho, F)$ is an essential subrepresentation of $(\pi,E)$ which
is isometric if $(\pi,E)$ is. Moreover, if $E$ is a reflexive Banach
space --- so that ${\cal B}(E)$ is a dual space --- and $\pi$ is
weak-weak$^\ast$ continuous, then so is $\rho$.
\end{lemma}
\begin{proof}
Straightforward.
\end{proof}
\par
For our next lemma, recall the notion of an 
ultrapower of a Banach space $E$ with respect to a (free)
ultrafilter $\cal U$ (\cite{Hei}); we denote it by $E_{\cal U}$.
\par
The lemma is a straightforward consequence of \cite[Proposition 5]{Daw}:
\begin{lemma} \label{matt}
Let $E$ be a superreflexive Banach space, and let $p \in
  (1,\infty)$. Then there is a free ultrafilter $\cal U$ such that the 
canonical representation of ${\cal B}(E)$ on $\ell_p(\posints, E)_{\cal U}$ is weak-weak$^\ast$ continuous.
\end{lemma}
\begin{lemma} \label{duals}
Let $G$ be a locally compact group, let $p,p' \in (1,\infty)$ be dual
to each other, and let $(\pi,E) \in \Rep_{p'}(G)$. Then, for each
$\phi \in \PF_{p',\pi}(G)$, there is a unique $g \in B_p(G)$ with $\|
g \|_{B_p(G)} \leq \| \phi \|$ such that
\begin{equation} \label{duality}
  \langle \pi(f), \phi \rangle = \int_G f(x) g(x) \, dx \qquad (f \in L_1(G)).
\end{equation}
Moreover, if $(\pi,E)$ is $p'$-universal, we have $\| g \|_{B_p(G)} = 
\| \phi \|$.
\end{lemma}
\begin{proof}
By Lemma \ref{matt}, there is a free ultrafilter such that the
canonical representation of $\PF_{p',\pi}(G)$ on
$\ell_{p'}(\posints,E)_{\cal U}$ is weak-weak$^\ast$ continuous. Use
Lemma \ref{subrep} to obtain an isometric, essential, and still
weak-weak$^\ast$ continuous subrepresentation $(\rho,F)$ of it.
\par
Since $E$ is a $\QSL_{p'}$-space and since the class of all
$\QSL_{p'}$-spaces is closed under the formation of $\ell_{p'}$-direct sums,
of ultrapowers, and of subspaces, $F$ is again a
$\QSL_{p'}$-space. Since $\rho$ is weak-weak$^\ast$ continuous and an
isometry, it follows that $\rho^\ast$ restricted to $F \ttensor_\pi
F^\ast$ is a quotient map onto $\PF_{p',\pi}(G)$. Let $\epsilon > 0$.
Then there are sequences $(\xi_n )_{n=1}^\infty$ in $F$ and
$(\psi_n)_{n=1}^\infty$ in $F^\ast$ such that
\[
  \| \phi \| \leq 
  \sum_{n=1}^\infty \| \xi_n \| \| \psi_n \| < \| \phi \| + \epsilon
  \qquad\text{and}\qquad
  \langle \rho(\pi(f)), \phi \rangle = \sum_{n=1}^\infty \langle \rho(f)\xi_n,
  \psi_n \rangle \quad (f \in L_1(G)).
\]
Since $\pi(L_1(G))$ is dense in $\PF_{p,\pi}(G)$, it follows that
$(\rho \circ \pi, F)$ is an essential representation of $L_1(G)$,
which therefore can be identified via (\ref{integral}) with an element
$(\sigma,F)$ of $\Rep_{p'}(G)$. Letting
\[
  g(x) := \sum_{n=1}^\infty \langle \sigma(x)\xi_n, \psi_n \rangle 
  \qquad (x \in G)
\]  
we obtain $g \in B_p(G)$ such that (\ref{duality}) holds. Moreover,
\[
  \| g \|_{B_p(G)} \leq \sum_{n=1}^\infty \| \xi_n \| \| \psi_n \| 
  < \| \phi \|+\epsilon;
\]
holds, and since $\epsilon > 0$ was arbitrary, this means that 
even $\| g \|_{B_p(G)} \leq \| \phi \|$.
\par
Suppose now that $(\pi,E)$ is $p'$-universal. Since the representation
of $L_1(G)$ induced by $(\pi,E)$ is essential, so is its infinite
amplification $(\pi^\infty, \ell_{p'}(\posints,E))$. With the
appropriate identifications in place, we thus have
\[  
  \ell_{p'}(\posints,E) \subset F \subset \ell_{p'}(\posints,E)_{\cal U}. 
\]
Consequently, $(\sigma, F)$ is also $p'$-universal. It then follows
from Lemma \ref{cyclem3} that $\| g \|_{B_p(G)} = \| \phi \|$.
\end{proof}
\par
In view of Lemma \ref{duals}, he following is now immediate:
\begin{theorem} \label{dualthm}
Let $G$ be a locally compact group,  and let $p,p' \in (1,\infty)$ be
dual to each other. Then:
\begin{items}
\item For any $(\pi,E) \in \Rep_{p'}(G)$, the dual space $\PF_{p',\pi}(G)^\ast$
  embeds contractively into $B_p(G)$.
\item The embedding of $\UPF_{p'}(G)^\ast$ into $B_p(G)$ is an isometric
  isomorphism.
\end{items}
\end{theorem}
\begin{remarks}
\item For $p=2$, the adverb ``contractively'' can be replaced by
``isometrically''. For $p \neq 2$, this is not true. To see this, assume
otherwise, and let $2 \leq q \leq p$ or $p \leq q \leq p$. 
Since $(\lambda_{q'}, L_{q'}(G)) \in \Rep_{p'}(G)$, we would thus have an
isometric embedding of $\PF_q(G)^\ast$ --- and thus of $A_q(G)$ --- into 
$B_p(G)$. For amenable $G$, this, in turn, would entail that $A_q(G) = A_p(G)$
holds isometrically. This is clearly impossible except in trivial cases.
\item As Michael Cowling pointed out to me, there is some overlap of this section with \cite{CF}. In particular, it is an immediate consequence of 
\cite[Theorem 2]{CF} that $B_p(G)$ is a dual Banach space.
\end{remarks}
\par
We conclude this section with a theorem that further clarifies
the relation between $B_p(G)$ and $A_p(G)$:
\begin{theorem}
Let $G$ be an amenable, locally compact group, and let $p,p' \in (1,\infty)$ be dual to each other. Then $\PF_{p'}(G)^\ast$, $B_p(G)$, and ${\cal M}(A_p(G))$ are equal with identical norms.
\end{theorem}
\begin{proof}
Since $G$ is amenable, we have $\PF_{p'}(G)^\ast = {\cal M}(A_p(G))$ with identical norms by \cite[Theorem 5]{Cow}, so that, by Theorem \ref{dualthm} and Corollary \ref{Fellcor2}, we have a chain
\[
  \PF_{p'}(G)^\ast \subset B_p(G) \subset {\cal M}(A_p(G)) =  \PF_{p'}(G)^\ast
\]
of contractive inclusions. This proves the claim.
\end{proof}
\begin{remark}
By \cite[Theorem 5]{Cow}, the equality $\PF_{p'}(G)^\ast = {\cal M}(A_p(G))$, even with merely equivalent and not necessarily identical norms, is also sufficient for the amenability of $G$. In view of the situation where $p =2$,
we suspect that $G$ is amenable if and only if $B_p(G) = {\cal M}(A_p(G))$ and if and only if $B_p(G) = \PF_{p'}(G)^\ast$.
\end{remark}
\dated
\vfill
\renewcommand{\baselinestretch}{1.2}
\begin{tabbing} 
{\it Address\/}: \= Department of Mathematical and Statistical Sciences \\
                 \> University of Alberta \\
                 \> Edmonton, Alberta \\
                 \> Canada, T6G 2G1 \\ \medskip
{\it E-mail\/}:  \> {\tt vrunde@ualberta.ca} \\[\medskipamount]
{\it URL\/}: \> {\tt http://www.math.ualberta.ca/$^\sim$runde/}  
\end{tabbing}


\begin{thebibliography}{L--N--R}
%
\begin{small}
%
\bibitem[Cow]{Cow} {\sc M.\ Cowling}, An application of Littlewood--Paley 
theory in harmonic analysis. {\it Math.\ Ann\/}.\  {\bf 241\/} (1979), 83--96.
%
\bibitem[C--F]{CF} {\sc M.\ Cowling} and {\sc G.\ Fendler}, On representations 
in Banach spaces. {\it Math.\ Ann.\/}\ {\bf 266\/} (1984), 307--315.  
%
\bibitem[Daw]{Daw} {\sc M.\ Daws}, Arens regularity of the algebra of
  operators on a Banach space. {\it Bull.\ London Math.\ Soc.\/}\ {\bf 36\/}
(2004), 493--503.
%
\bibitem[D--F]{DF} {\sc A.\ Defant} and {\sc K.\ Floret}, {\it Tensor Norms and Operator Ideals\/}. North-Holland, 1993. 
%
\bibitem[E--R]{ER} {\sc E.\ G.\ Effros} and {\sc Z.-J.\ Ruan}, {\it Operator Spaces\/}. Clarendon Press, Oxford, 2000.
%
\bibitem[Eym]{Eym} {\sc P.\ Eymard}, L'alg\`ebre de Fourier d'un groupe localement compact. {\it Bull.\ Soc.\ Math.\ France\/} {\bf 92\/} (1964), 181--236.   
%
\bibitem[For 1]{For1} {\sc B.\ E.\ Forrest}, Arens regularity and the $A_p(G)$ 
algebras. {\it Proc.\ Amer\. Math.\ Soc.\/}\ {\bf 119\/} (1993), 595--598.
%
\bibitem[For 2]{For2} {\sc B.\ E.\ Forrest}, Amenability and the structure of 
the algebras $A_p(G)$. {\it Trans.\ Amer.\ Math.\ Soc.\/}\ {\bf 343\/} (1994),  233--243
%
\bibitem[G--L]{GL} {\sc I.\ Glicksberg} and {\sc K.\ de Leeuw}, The decomposition of certain group representations. {\it J.\ Anal.\ Math.\/}\ {\bf 15\/} (1965), 135--192.
%
\bibitem[Hei]{Hei} {\sc S.\ Heinrich}, Ultraproducts in Banach space theory. {\it J.\ reine angew.\ Math.\/}\ {\bf 313\/} (1980), 72--104
%
\bibitem[Her 1]{Her1} {\sc C.\ Herz}, The theory of $p$-spaces with an
  application to convolution operators. {\it Trans.\ Amer.\ Math.\ Soc.\/}\
  {\bf 154\/} (1971), 69--82.
%
\bibitem[Her 2]{Her2} {\sc C.\ Herz}, Harmonic synthesis for subgroups. {\it Ann.\ Inst.\ Fourier\/} (Grenoble) {\bf 23\/} (1973), 91--123.
%
\bibitem[H--R]{HR} {\sc C.\ Herz} and {\sc N.\ Rivi\`ere}, Estimates for translation-invariant operators on spaces with mixed norms.
{\it Studia Math.\/}\ {\bf 44\/} (1972), 511--515.
%
\bibitem[J--M]{JM} {\sc P.\ Jaming} and {\sc W.\ Moran}, Tensor products and $p$-induction of representations on Banach spaces. {\it Collect.\ Math.\/}\ {\bf 51\/} (2000), 83--109.
%
\bibitem[Kwa]{Kwa} {\sc S.\ Kwapie\'n}, On operators factoring through $L_p$-space. {\it Bull.\ Soc.\ Math.\ France, M\'em.\/}\ {\bf 31--32\/} (1972), 215--225.
%
\bibitem[LeM]{LeM} {\sc C.\ LeMerdy}, Factorization of $p$-completely bounded
  multilinear maps. {\it Pacific J.\ Math.\/}\ {\bf 172\/} (1996), 187--213.
%
\bibitem[L--N--R]{LNR} {\sc A.\ Lambert}, {\sc M.\ Neufang}, and {\sc V.\
  runde}, Operator space structure and amenability for Fig\`a-Talamanca--Herz
  algebras. {\it J.\ Funct.\ Anal.\/}\ {\bf 211\/} (2004), 245--269.
%
\bibitem[L--R]{LR} {\sc J.\ Lindenstrauss} and {\sc H.\ P.\ Rosenthal}, The
  ${\cal L}_p$ spaces. {\it Israel J.\ Math.\/}\ {\bf 7\/} (1969), 325--349.
%
\bibitem[Loh]{Loh} {\sc N.\ Lohou\'e}, Estimations $L^p$ des coefficients de 
repr\'esentation et op\'erateurs de convolution. {\it Adv.\ in Math.\/}\ 
{\bf 38\/} (1980), 178--221.
%
\bibitem[Mia]{Mia} {\sc T.\ Miao}, Compactness of a locally compact group $G$ 
and geometric properties of $A_p(G)$. {\it Canad.\ J.\ Math.\/}\ {\bf 48\/} 
(1996), 1273--1285.
%
\bibitem[Pie]{Pie} {\sc J.\ P.\ Pier}, {\it Amenable Locally Compact Groups\/}. Wiley-Interscience, 1984. 
%
\end{small}
%
\end{thebibliography}
\end{document}